\documentclass[lettersize,journal]{IEEEtran}
\usepackage{amsmath,amssymb,amsfonts}
\usepackage{algorithmic}
\usepackage{array}
\usepackage[caption=false,font=normalsize,labelfont=sf,textfont=sf]{subfig}
\usepackage{textcomp}
\usepackage{filecontents}
\usepackage{stfloats}
\usepackage{url}
\usepackage{booktabs}
\usepackage{comment}
\usepackage{verbatim}
\usepackage{graphicx}
\usepackage{cite}

\usepackage[ruled,vlined]{algorithm2e}
\makeatletter
\let\NAT@parse\undefined
\makeatother
\usepackage[colorlinks,allcolors=blue]{hyperref}
\makeatletter
\renewcommand*{\eqref}[1]{%
  \hyperref[{#1}]{\textup{\tagform@{\ref*{#1}}}}%
}
\makeatother

\def\BibTeX{{\rm B\kern-.05em{\sc i\kern-.025em b}\kern-.08em
    T\kern-.1667em\lower.7ex\hbox{E}\kern-.125emX}}

\begin{document}

\title{{\fontsize{24}{30}\selectfont A Lyapunov-based MPC for Distributed Multi Agent Systems with Time Delays and Packet Dropouts using Hidden Markov Models}}

\author{Loaie Solyman$^{1}$, Aamir Ahmad$^2$, Ayman El-Badawy$^{1,*}$ 
\thanks{This work has been supported by the "Optimal Networked Control" Project funded by the German Academic Exchange Service (DAAD) and the Federal Ministry for Education, Research, Technology and Aerospace (BMFTR) in the framework of the DAAD-TNB-BINA funding project "GUC: Building A Sustainable Future" (project ID: 57710434).}
\thanks{$^{1}$Mechatronics Engineering Department, Faculty of Engineering and Materials Science, German University in Cairo, Cairo, Egypt}%
\thanks{$^{2}$Institute of Flight Mechanics and Control (IFR), University of Stuttgart, Stuttgart, Germany}
\thanks{$^*$Corresponding author (ayman.elbadawy@guc.edu.eg)}}

\markboth{IEEE Transactions on Control Systems Technology, VOL. XX, NO. XX, XXXX 2017}
{Solyman \MakeLowercase{\textit{et al.}}: A LMPC for Distributed MASs with Time Delays and Packet Dropouts using HMMs}

\maketitle

\begin{abstract}
We propose a SCHMM-LMPC framework, integrating Semi-Continuous Hidden Markov Models with Lyapunov-based Model Predictive Control, for distributed optimal control of multi agent systems under network imperfections. The SCHMM captures the stochastic network behavior in real time, while LMPC ensures consensus and optimality via Linear Matrix Inequalities (LMIs). The developed optimal control problem simultaneously minimizes three elements. First, the control effort is reduced to avoid aggressive inputs and second, the network-induced error caused by time delays and packet dropouts. Third, the topology-induced error, as the distributed graph restricts agents’ access to global information. This error is inherent to the communication graph and cannot be addressed through offline learning. To overcome this, the study also introduces the incremental Expectation Maximization (EM) algorithm, enabling online learning of the SCHMM. This adaptation allows the framework to mitigate both network and topology errors while maintaining optimality through MPC. Simulations validate the effectiveness of the proposed SCHMM-LMPC, demonstrating adaptability in multi agent systems with diverse topologies. 
\end{abstract}

\section{Introduction}

\IEEEPARstart{M}{ulti} Agent Systems (MASs) have gained considerable attention in fields such as robotics, autonomous vehicles, smart grids, and distributed sensor networks \cite{4n}. These systems rely on consensus-oriented control, where agents coordinate their actions through local connection \cite{5n}. Achieving consensus, however, is difficult under constraints like time delays, packet dropouts \cite{36}, and distributed topologies \cite{20}. Such issues are critical in team formations, where delays or data loss hinder coordination, or in smart grids, where distributed agents must operate synchronously under diverse constraints. Despite extensive research, no unified solution fully addresses these challenges \cite{38}.

The study of distributed optimal control in MASs remains challenging due to three interrelated issues, namely, distributed topologies \cite{40}, the need for optimality \cite{41}, and network imperfections \cite{40}. In distributed settings, agents lack centralized control, exchanging only local data with neighbors \cite{40}. This improves scalability and resilience but complicates control law design, as agents must act on incomplete global information. The optimality issue requires control laws to meet performance criteria such as constrained control effort \cite{45}. Achieving this collectively, with agents solving coupled optimization problems under partial information, necessitates distributed optimization algorithms \cite{44}. Network imperfections, such as delays and packet losses, further impair information exchange, degrading performance or causing instability \cite{46}. These stochastic effects are hard to model, requiring probabilistic strategies. The three issues are tightly coupled, as unreliable communication complicates distributed decision-making, while ensuring optimality demands precision often disrupted by network imperfections \cite{47}. Thus, comprehensive solutions must address all three jointly. This paper develops such a framework, systematically incorporating distributed structure, guaranteeing optimality, and explicitly accounting for delays and packet dropouts.

Regarding topology, three paradigms exist, the centralized, decentralized, and distributed. Centralized schemes rely on a single controller, which may achieve global optimality but suffer from scalability and single-point failure \cite{48}. Decentralized schemes let agents act independently with only local data, enhancing resilience but yielding suboptimal outcomes \cite{49}. Distributed topologies strike a balance where agents interact with neighbors according to a network graph and compute local actions using both local and shared data \cite{42}. This approach offers scalability and resilience to single-point failure, making it the most practical choice for MASs in uncertain environments \cite{50}.

The optimality issue extends beyond mere consensus, requiring control strategies that minimize a defined cost while respecting constraints. In dynamic or resource-limited environments, suboptimal solutions reduce efficiency. Model Predictive Control (MPC) addresses this by optimizing a cost function over a prediction horizon under constraints \cite{51}. Coupled with Lyapunov-based analysis, MPC guarantees optimal performance and consensus, making it effective for MAS applications \cite{52}.

Network imperfections such as time delays and packet dropouts present further challenges. Delays hinder timely decisions, while dropouts cause data loss \cite{21}. Mitigation strategies include delay-compensation methods, predictive control, and event-triggered schemes \cite{20}. More recent works embed stochastic models like Hidden Markov Models (HMMs) into consensus algorithms, providing probabilistic guarantees \cite{36}. Study \cite{36} introduced a Semi-Continuous HMM (SCHMM) called the Single Model Scheme (SMS) to jointly capture delays and dropouts. However, this work only addressed single control systems and not MASs.

Some studies attempt to combine subsets of the three problems. For example, \cite{53} addressed distributed topology and optimality via a randomized Jacobi proximal Alternating Direction Method of Multipliers (ADMM), but the method is complex and ignores communication constraints. Another study tackled topology and communication jointly \cite{54}, requiring agents to track all possible packet loss scenarios. While effective, the assumptions are local to limited geographical areas as mentioned by the study itself. Other contributions like \cite{55} and \cite{56} similarly address only two of the three problems. Where \cite{55} used a Markov chain to predict the behavior of packet dropouts only, solving the distributed topology issue along with partial solution to the networks issue. Also, \cite{56} utilized iterative learning to predict the behavior of continuous data losses while adhering to a cost function, solving only the optimality issue and partially attacking the networks issue. No existing study comprehensively solves distributed optimal control in MASs under both time delays and packet dropouts. These mentioned studies also have a drawback, which is their formulation of the problem. The distributed nature of the topology constraints the number of neighbors each agent can communicate with, in addition, this communication is impaired by the network. Earlier studies formulate the problem in a way that heavily couples the topological and network issues, which then poses a difficult problem to solve. This compels them to make conservative assumptions on the topology or the network. This drawback was noticed in the literature and hence this study would overcome that through an independent formulation of the problem.

This paper introduces the Semi-Continuous Hidden Markov Model – Lyapunov-based Model Predictive Control (SCHMM-LMPC), the first framework to simultaneously address distributed topology, optimality, and network uncertainties. This method groups neighboring agents into subsystems which reduces the computational burden. This extracts a compact system, where an error term linked to distributed topology, is introduced. Unlike prior work where SCHMMs only modeled network errors \cite{36}, this study considers errors also arising from topology. To handle this, an incremental Expectation Maximization (EM) algorithm that is novel in SCHMM contexts, is introduced for online learning. Moreover, LMPC ensures optimality and consensus through Lyapunov-based guarantees.

The problem at hand is non-trivial because it deals with two apparently dependent error terms. The first is the error due to the network delays and losses. And, the second is the error due to the distributed topology, as each agent can only communicate with a subset of the graph. The main contribution of this study is a compact system formulation that removes this interdependence. By leveraging the SCHMM prediction mechanism to handle network effects first then grouping each agent with its immediate neighbors into a local subsystem. This scheme tackles the error from the network by predicting the states of the neighbors beforehand, and tackles the error due to the distributed topology by using those predictions to drive each agent to its local consensus point. This scheme would converge all local consensus points to the global consensus point in a finite time. So, this new formulation treats the two error terms independently and this, in turn, enables more relaxed and practical constraints on delays and topologies compared to previous approaches.


The main contributions are summarized in two points. The first is the compact system formulation which reduces computation while capturing the distributed error term, a Lyapunov-based MPC ensures optimality and consensus via Linear Matrix Inequalities (LMIs). The second is the use of SCHMM to predict network behavior and compensate for delays and dropouts, extended with incremental EM to also capture topology-induced errors. This marks the first integration of HMMs into MPC-based MAS solutions. Effectiveness is validated through numerical examples, showing consensus under adverse conditions and outperforming methods addressing only subsets of the challenges.

The paper is organized as follows. Section \ref{sec1} formulates the problem and challenges. Section \ref{sec2} describes the methodology namely the compact system formation, LMPC, and SCHMM integration. Section \ref{sec3} presents theoretical analysis. Section \ref{sec4} provides simulations, and Section \ref{sec5} concludes with findings and contributions.

\subsection{Preliminaries}

Standard notation is used throughout this paper. Let $\mathbb{R}$ and $\mathbb{N}$ denote the sets of real and natural numbers (excluding zero), respectively. $\mathbb{R}^n$ is the $n$-dimensional Euclidean space. $1_n$ and $0_n$ denote column vectors of ones and zeros. For a matrix $A$, $A^T$ and $A^{-1}$ represent its transpose and inverse. $P>0$ implies that $P$ is real, symmetric, and positive definite. $I$ and $\textbf{0}$ denote the identity and zero matrices of appropriate dimensions. The symbol $*$ denotes symmetric blocks in partitioned matrices. A Schur stable matrix has eigenvalues $\lambda_i$, $i=\{1,2,\dots,n\}$, lying within the unit circle. $P(A)$ denotes the probability of an event $A$, $\otimes$ the Kronecker product, $A*B$ the element-wise matrix multiplication, and $(A)_{i*}$ the $i$th row of $A$. The set cardinality is denoted $\mathbf{card}(A)$, and the Hadamard inverse of $A$ (with $a_{ij}>0$) is $A^{\circ(-1)}=(1/a_{ij})$.

Consider a MAS with $N$ agents whose communication is represented by a graph $\mathcal{G}(\mathcal{V}, \mathcal{E}, \mathcal{A})$ of order $N$, where $\mathcal{V} = \{v_1, v_2, \dots, v_N\}$, $\mathcal{E} = \{e_{ij}=(v_i,v_j)\} \subset \mathcal{V}\times\mathcal{V}$, and $\mathcal{A}=[a_{ij}] \in \mathbb{R}^{N \times N}$. Here, $a_{ij}=1$ if $(v_j,v_i)\in\mathcal{E}$, and $a_{ij}=0$ otherwise, with $a_{ii}=0$. Each edge $e_{ij}$ indicates that agent $j$ receives information from agent $i$. The neighbor set of node $v_i$ is $\mathcal{N}_i = \{v_j \mid e_{ji}=(v_j,v_i)\in\mathcal{E}\}$. An undirected graph contains a spanning tree if there exists a root node connected to all others through undirected paths. The Laplacian matrix is $L = D - \mathcal{A}$, where $D = \mathrm{diag}\{\sum_{j=1}^N a_{1j}, \sum_{j=1}^N a_{2j}, \dots, \sum_{j=1}^N a_{Nj}\}$.

\section{Problem Formulation}
\label{sec1}
Consider a MAS with $N$ heterogeneous agents where each agent $i \in \{1,2,3,\hdots,N\}$ can have the following general linear dynamics:

\begin{figure*}[!t]
\centering
\includegraphics[scale=0.35]{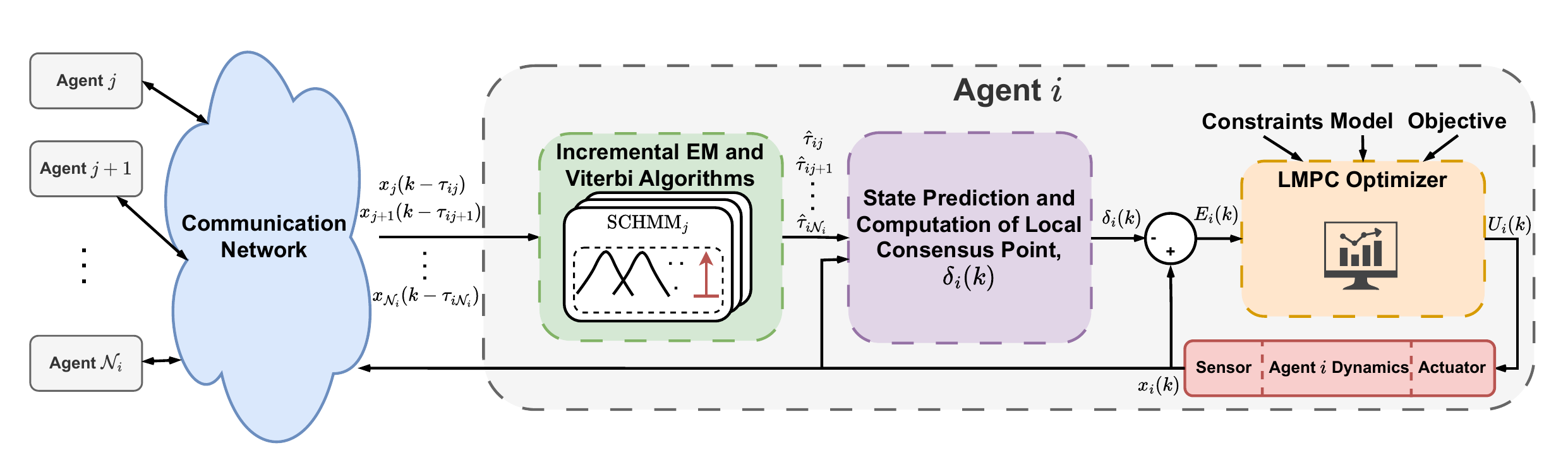}
\caption{SCHMM-LMPC framework for agent $i$ describing the data flow of the received delayed information to drive agent $i$ towards global consensus.}
\label{ffig}
\end{figure*}

\begin{equation} \label{eq1}
x_i(k+1)=A_ix_i(k)+B_iu_i(k)
\end{equation}

\noindent where $x_i(k) \in \mathbb{R}^n$ and $u_i(k) \in \mathbb{R}^m$ denote the state and input of agent $i$, respectively. And, $A_i \in \mathbb{R}^{n \times n}$ and $B_i \in \mathbb{R}^{n \times m}$ denote the system and input matrices of agent $i$, respectively. And, the consensus of the agents is achieved in this study by acknowledging that translational states, $\{x_T\}$, of all agents tend asymptotically to convergence. However, the distributed nature of the topology dictates that some agents might not be in connection with all other agents. And thus, the problem is not reduced to achieving consensus on an agreed-upon-point in space as not all agents have access to global information. Nevertheless, all graphs used in this study are assumed to be undirected as well as containing a spanning tree.

The SCHMM-LMPC framework can be inspected in Fig. \ref{ffig}, where it is assumed that the information received from neighbors of agent $i$ is delayed by $\tau_{ij}$. And $\tau_{ij}$ denotes the overall time delay imposed on the packet of information sent from the sending agent until it reaches the destination agent. The word 'overall' used in the previous sentence is to consider the effect of time delays and packet dropouts together. It has to be noted as well, that agent $i$ broadcasts its information to its neighbors through the communication network and thus the framework implemented at agent $i$ shown in Fig. \ref{ffig} is also implemented for all agents in the MAS under study.

Now, in this study, the use of the SCHMM will facilitate the computation of the predicted states of the neighboring agents as will be discussed later in this study. This can be examined in Fig. \ref{ffig} where the novel incremental EM algorithm for online training of the SCHMM is used. It has to be noted here that there exists, in the framework, an SCHMM for each neighbor, as interactions with neighbors differ through the simulation giving rise to a different set of delays for each neighbor. This dictates that the incremental EM and Viterbi algorithms take in the delayed information from the neighbors and output the predicted delays for each neighbor. This then follows that $x_j(k-\tau_{ij})$, the received delayed information, along with $\hat{\tau}_{ij}$, the predicted delay experienced by agent $i$ upon receiving a packet from agent $j$, can be used to arrive at $\hat{x}_j(k)$, the predicted instantaneous position of agent $j$, for all neighboring agents $j$. This then goes to help in the computation of the instantaneous state prediction of the neighbors, as well as the local consensus point, $\delta_i$, which is the instantaneous target location for agent $i$. 

With this information at hand, the framework shown in Fig. \ref{ffig} presents the use of the LMPC which is responsible for the computation of the optimal control input for agent $i$, $U_i(k)$. It takes as input the error computed from the position of agent $i$ and $\delta_i$. Governed by a Lyapunov candidate function, Fig. \ref{ffig} shows the LMPC which complies with an objective, simulates the model dynamics and maintains adherence to constraints in order to arrive at the optimal input commands, $U_i(k)$. The specific details of this optimization problem will follow in this study. The optimal input commands are then actuated upon in the physical plant of agent $i$ which then sends its information to its neighbors, as shown in Fig. \ref{ffig}.

\section{Methodology}
\label{sec2}

We first address the distributed topology by grouping each agent with its immediate neighbors into local subgroups and deriving the corresponding compact system to establish consensus convergence. We then consider the second challenge, that is the unknown, time-varying delays and packet dropouts, which is handled through the SCHMM and its role in mitigating network-induced imperfections.


\subsection{Grouping}

It has been discussed earlier in this study that for a MAS with $N$ agents, the set $\mathcal{N}$ be the set containing all agents in our MAS and thus $\mathbf{card}(\mathcal{N}) = N$. Moreover, let us define $N$ new subsets, $\mathcal{N}_i \hspace{2mm} \forall ~ i \in \mathcal{N}$. These subsets include each agent and its immediate neighbors. 

\begin{figure}[!t]
\centering
\includegraphics[scale=0.3]{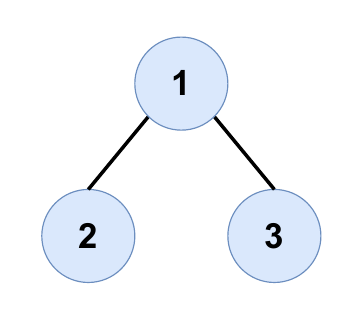}
\caption{Example Multi Agent System with three agents.}
\label{fig1}
\end{figure}

\textbf{Example:} You may refer to Fig. \ref{fig1} and Table \ref{tab1} for clarification of the sets $\mathcal{N}$ and $\mathcal{N}_i$ along with their corresponding dynamics matrices $A_c$ and $B_c$ which are used in later mathematical formulation. These $N$ new sets will be the cornerstone of tackling the problem of the distributed topology as now the graph has been divided into $N$ graphs that do not exhibit the distributed topology issue which will help pave the way towards a solution for this problem. This goes to show that this solution is easily applicable to MASs with any topology, specially the distributed topology at hand.

In order to facilitate distributed consensus, we define a new term $\delta_i(k)$ for each agent $i$ as follows:

\begin{equation} \label{eq81}
\delta_i(k)=\frac{1}{\textbf{card}(\mathcal{N}_i)}(x_i(k)+\sum_{j=1}^{\mathcal{N}_i}\hat{x}_j(k))
\end{equation}

\begin{table}[!t]
    \centering
\caption{Sets $\mathcal{N}$ and $\mathcal{N}_i$, and matrices $A_c$ and $B_c$, for the agents in the MAS shown in Fig. \ref{fig1}.}
    \begin{tabular}{@{}llll@{}}
\toprule
MAS/Agent & Set&  $A_c$ & $B_c$\\
\midrule
MAS in \\ Fig. \ref{fig1}  & $\mathcal{N} \hspace{1mm} =\{1,2,3\}$ & - & - \\
1    & $\mathcal{N}_1=\{1,2,3\}$ & $\text{diag}\{A_1,A_2,A_3\}$ &  $\text{diag}\{B_1,B_2,B_3\}$\\
2   & $\mathcal{N}_2=\{1,2\}$ & $\text{diag}\{A_1,A_2,0\}$ & $\text{diag}\{B_1,B_2,0\}$\\
3   & $\mathcal{N}_3=\{1,3\}$ & $\text{diag}\{A_1,0,A_3\}$  & $\text{diag}\{B_1,0,B_3\}$\\
\bottomrule
\end{tabular}
    
    \label{tab1}
\end{table}

\noindent where the term $\delta_i(k)$ denotes the local consensus point. This is defined as the target position for each agent $i$ computed from the position of agent $i$ along with its immediate neighbors $j$ with $j \neq i$. The term $\delta_i(k)$ will act as the target position of agent $i$ that will drive all agents into the state of consensus. However, as our study focuses as well on the network imperfections, the positions of the neighbors should have been the delayed form of the position, $x_j(k-\tau_{ij})$. But, this study mitigates the effects of time delays and packet dropouts though the use of HMMs, which will be explained later in this study. This facilitates agent $i$ to have a prediction of each neighbor's overall delay, and thus computes a prediction of the instantaneous position of the neighbors, $\hat{x}_j(k)$ which is then used to calculate $\delta_i(k)$ as shown in \eqref{eq81}.

The advantage behind defining such a term is seen through the fact that each agent $i$ is now solving their version of the global problem. And by virtue of having a spanning tree in all graphs studied in this paper, the global consensus is guaranteed when all agents converge to their respected $\delta_i(k)$.

\begin{figure}[!t]
\centering
\includegraphics[width=0.3\textwidth,height=0.2\textwidth]{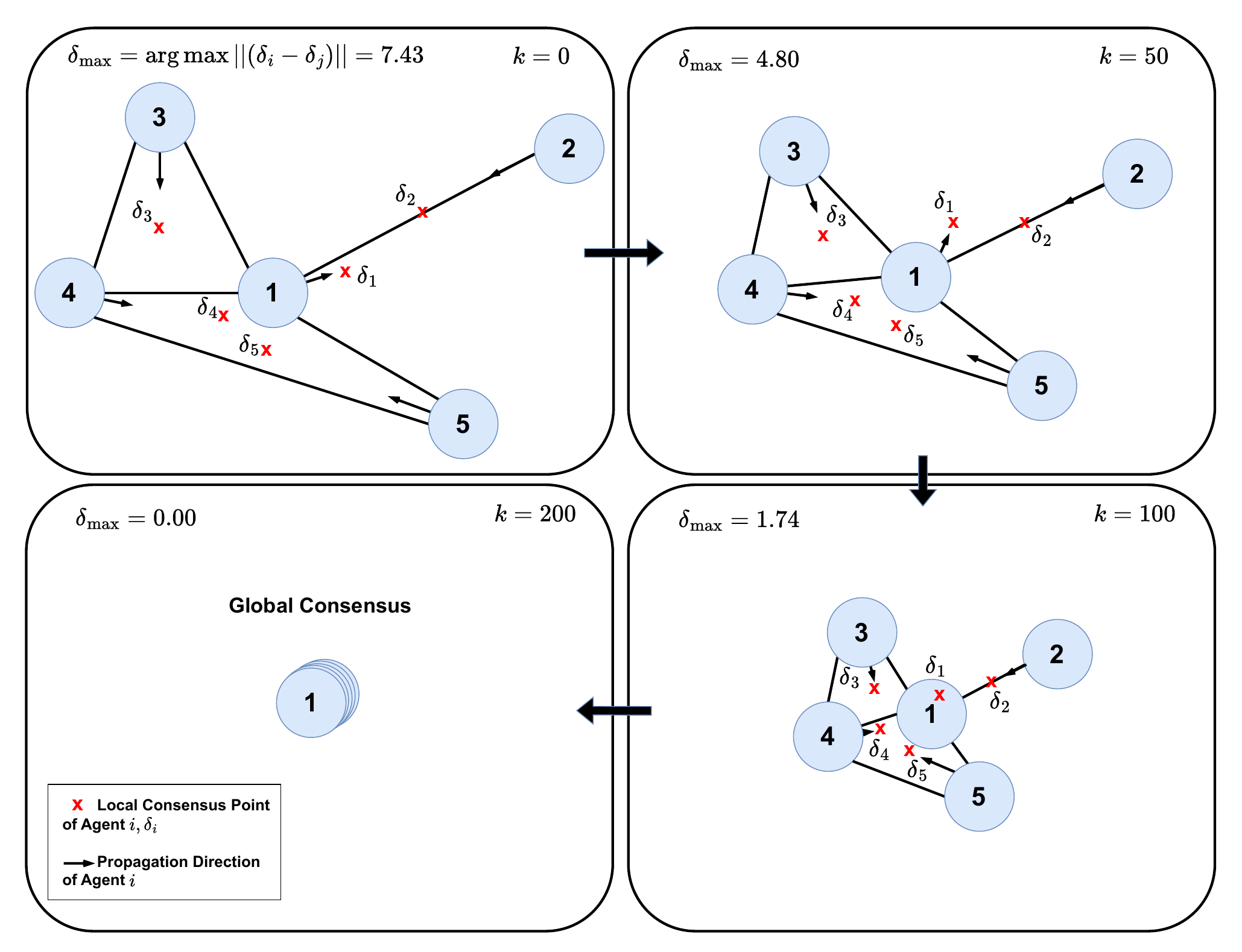}
\caption{Example Multi Agent System with five agents (starting top-left) showing local consensus points to scale computed using \eqref{eq81} and propagation directions of the agents achieving global consensus.}
\label{ffig2}
\end{figure}

\textbf{Example:} It is advised to refer to Fig. \ref{ffig2} which further elaborates the use of the local consensus point, $\delta_i(k)$. Fig. \ref{ffig2} presents the example of a MAS with five agents highlighting two key elements to aid in visualizing the evolution of the MAS applying the local consensus point mechanism. The first element is the local consensus point of each agent, $\delta_i$, which is the instantaneous target location for each agent. The second element is the propagation direction of each agent denoting the direction in which the agent actuates, which is towards the local consensus point. The term $\delta_{max}$ is defined as the maximum difference between any two local consensus points in the MAS. Fig. \ref{ffig2} shows the evolution of $\delta_{max}$ converging to zero, i-e the MAS achieving global consensus, even though it relies on local consensus points.

Now, that the term $\delta_i(k)$ has been defined, it is wise to now define the error to be:

\begin{equation} \label{eq1000}
e_i(k)=x_i(k)-\delta_i(k)
\end{equation}

As of now, the consensus error for each agent $i$ has been defined using the predicted states from the SCHMM, and the problem could now be formulated as, $\lim_{k\rightarrow \infty}e_i=0 ~ \forall ~ i \in \mathcal{N}$. And this can be used to define the cost function for each agent $i$ using the error term in \eqref{eq1000} as:

\begin{equation} \label{eq1001}
j_i(k)=\sum_{p=1}^{N_{\alpha}}e_i^T(k+p)Pe_i(k+p)+u_i^T(k+p)Qu_i(k+p)
\end{equation}

\noindent where $0<P \in \mathbb{R}^{n \times n}$ and $0<Q \in \mathbb{R}^{m \times m}$ denote positive definite matrices that need to be appropriately selected. Also, $N_{\alpha}$ is the optimization horizon, which varies from time step to the other according to a constant $\alpha > 0$ that is defined as the upper bound for the chosen Lyapunov candidate function which will be chosen later in this study. Consequently, minimizing the cost function $j_i$ for each agent $i$ ensures the fulfillment of the previously defined problem, thereby guaranteeing consensus among all agents.

\subsection{Compact Form}

Now, that the groups have been defined, our method of control design needs to be considered. This will be achieved through the augmentation of all agents into a compact form which will then be used further in the Lyapunov based MPC.

Using \eqref{eq1} and augmenting the system for all $N$ agents, one gets:

\begin{equation} \label{eq8}
X(k+1)=A_mX(k)+B_mU(k)
\end{equation}

\noindent which is defined using the following:

\begin{equation} \label{eq9}
\begin{gathered}
X(k)=[x_1^T(k), x_2^T(k), x_3^T(k), \hdots, x_N^T(k)]^T,\\
U(k)=[u_1^T(k), u_2^T(k), u_3^T(k), \hdots, u_N^T(k)]^T,\\
A_m=\text{diag}\{A_1,A_2,A_3,\hdots,A_N\}, \\B_m=\text{diag}\{B_1,B_2,B_3,\hdots,B_N\}
\end{gathered}
\end{equation}

\noindent given that $X(k) \in \mathbb{R}^{nN}$ and $U(k) \in \mathbb{R}^{mN}$ denoting the states and inputs of the global system of the MAS. With $A_m \in \mathbb{R}^{nN \times nN}$ and $B_m \in \mathbb{R}^{nN \times mN}$ denoting the system and input matrices of the global system. 

However, to address the distributed topology nature of the MAS, we need to define local versions of the global problem at each agent $i$ as follows:

\begin{equation} \label{eq8}
\begin{gathered}
X_i(k+1)=(I_N*(1_N^T \otimes (\mathcal{A}+I_N)_{i*}^T) \otimes I_N)A_mX_i(k)+\\(I_N*(1_N^T \otimes (\mathcal{A}+I_N)_{i*}^T) \otimes I_N)B_mU_i(k)
\end{gathered}
\end{equation}

\noindent which is defined using the following:

\begin{equation} \label{eq9}
\begin{gathered}
X_i(k)=[\hat{x}_1^T(k), \hdots, \hat{x}_{i-1}^T(k), x_i^T(k), \hat{x}_{i+1}^T(k), \hdots, \hat{x}_N^T(k)]^T,\\
U_i(k)=[u_1^T(k), \hdots, u_{i-1}^T(k), u_i^T(k), u_{i+1}^T(k), \hdots, u_N^T(k)]^T
\end{gathered}
\end{equation}

\noindent given that $X_i(k) \in \mathbb{R}^{nN}$ and $U_i(k) \in \mathbb{R}^{mN}$ denoting, respectively, the predicted states and inputs of the global system from the point of view of agent $i$ of the MAS along with its own local data. With the operator $(I_N*(1_N^T \otimes (\mathcal{A}+I_N)_{i*}^T) \otimes I_N)$ indicating the extraction of information of the neighbors of agent $i$. To simplify the notation, let $A_c=(I_N*(1_N^T \otimes (\mathcal{A}+I_N)_{i*}^T) \otimes I_N)A_m$ and $B_c=(I_N*(1_N^T \otimes (\mathcal{A}+I_N)_{i*}^T) \otimes I_N)B_m$ to arrive at a simpler notation as follows:

\begin{equation} \label{eq10}
X_i(k+1)=A_cX_i(k)+B_cU_i(k)
\end{equation}

And to make the notation clearer, the example shown in Fig. \ref{fig1} and Table \ref{tab1} describe the matrices $A_c$ and $B_c$ for a simple MAS. The distributed nature of the topology is still maintained in \eqref{eq10} through the adjacency matrix, $\mathcal{A}$ which describes the connection in the graph. However, this compact form did offer an advantage of having the term $X_i(k)$ for all agents in the MAS, to be of equal vector size which helps in the mathematical manipulation further.

Now, using \eqref{eq81} and \eqref{eq1000}, we consider the compact error expression for each agent $i$ as follows:

\begin{equation} \label{eq12}
\begin{gathered}
    E_i(k)=((-(A_p1_{1 \times N})^{\circ(-1)}A_q) \otimes I_n)X_i(k)
\end{gathered}
\end{equation}

\noindent with $A_p=[
        \textbf{card}(\mathcal{N}_1),
        \textbf{card}(\mathcal{N}_2),
        \hdots,
        \textbf{card}(\mathcal{N}_N)]^T$ and $A_q=\text{diag}\{(1-\textbf{card}(\mathcal{N}_1)),(1-\textbf{card}(\mathcal{N}_2)), \hdots, (1-\textbf{card}(\mathcal{N}_N))\}$, and letting $A_e=((-(A_p1_{1 \times N})^{\circ(-1)}A_q) \otimes I_n)$. It has to be noted here that the expression $(-(A_p1_{1 \times N})^{\circ(-1)}A_q)$ yields a full rank matrix that is multiplied by $I_n$ in a kronecker product. And as the first operator of the product is a full rank invertible matrix, then the output of the product, $A_e$, is full rank invertible matrix as well by virtue of $\det(A_{n \times n} \otimes I_{n \times n})=\det(A_{n \times n})^n \det(I_{n \times n})^n=\det(A_{n \times n})^n$. Then the compact error expression for each agent $i$ is simplified to: 

\begin{equation} \label{eq13}
E_i(k)=A_eX_i(k)
\end{equation}

It has to be noted here that if the system was not impaired by the network imperfections at all, then this error still pertains to the system. Then, it has to be declared that this error term shown in \eqref{eq13} is an error due to the distributed nature of the system.

Now, we solve for the $(k+1)^{\text{th}}$ time instant and using \eqref{eq10} and \eqref{eq13} to get:

\begin{equation} \label{eq14}
\begin{gathered}
E_i(k+1)=A_eX_i(k+1)=A_eA_cX_i(k)+A_eB_cU_i(k)
\end{gathered}
\end{equation}

\noindent and now from \eqref{eq13} and by the invertibility property stated earlier, we also get:

\begin{equation} \label{eq15}
X_i(k)=A_e^{-1}E_i(k)
\end{equation}

\noindent now using \eqref{eq14} and \eqref{eq15} we get:

\begin{equation} \label{eq16}
\begin{gathered}
E_i(k+1)=A_eA_cA_e^{-1}E_i(k)+A_eB_cU_i(k)
\end{gathered}
\end{equation}

The input command expression $U_i(k)$ will now be formulated using MPC-based state feedback. Consequently, the input command $U_i(k)$ will ensure that the states of the agents converge, or equivalently, that $E_i(k)$ approaches zero for all agents $i$. 

The corresponding compact form of the cost function is given as follows:

\begin{equation} \label{eq20}
J_i(k)=\sum_{p=0}^{N_{\alpha}}E_i^T(k+p)P_JE_i(k+p)+U_i^T(k+p)Q_JU_i(k+p)
\end{equation}

\noindent with $P_J=\text{diag}\{P,P,P,\hdots,P\}$ and $Q_J=\text{diag}\{Q,Q,Q,\hdots,Q\}$ and $J_i$ is to be minimized. Then we define the input command as:

\begin{equation} \label{eq21}
U_i(k)=KE_i(k)
\end{equation}

\noindent and $K$ is to be designed using Lyapunov based MPC. This then completes the compact form of the system to tackle the problem of the distributed nature of the topology of the MAS under study.

\subsection{Network Imperfections}

The task of network modeling consists of developing a model, $\lambda$, capable of replicating the actual network’s behavior. This challenge is divided into two subproblems. The first involves gathering raw time delay and packet dropout data via simulations, which can be accomplished using NS-2 which is the most widely utilized network simulator in the literature \cite{3}. The second involves adjusting the SCHMM through the EM algorithm, followed by executing the model with the Viterbi algorithm to forecast new time delays and packet losses consistent with the training set employed during tuning.

\begin{figure}[!t]
\centering
\includegraphics[trim={0 0 5.5cm 0},clip,scale=0.2]{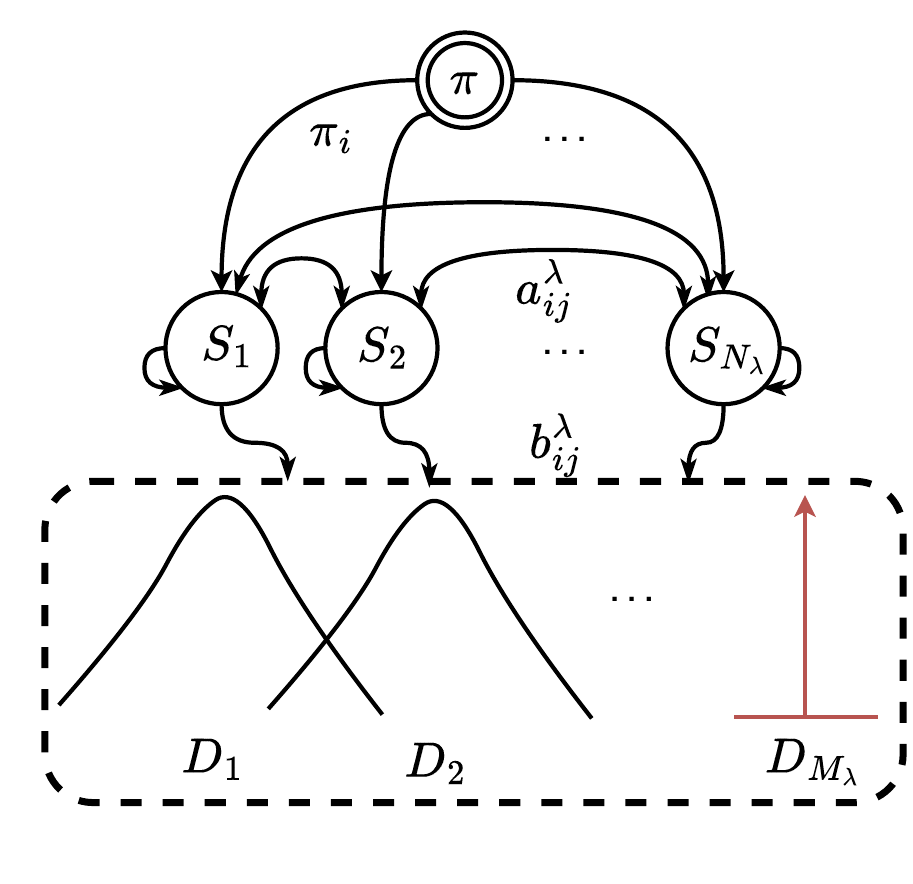}
\caption{Internal structure of the Semi Continuous Hidden Markov Model.}
\label{fig6}
\end{figure}

The SCHMM, illustrated in Fig. \ref{fig6}, can be characterized as a standard HMM \cite{1}. It can be represented using the following tuple:

\begin{equation} \label{eq6}
\lambda = (N_{\lambda},M_{\lambda},A,B,\pi)
\end{equation}

\noindent where, $N_{\lambda} \in \mathbb{N}$ represents the number of hidden or latent states, depicted as circles in Fig. \ref{fig6}, while $M_{\lambda} \in \mathbb{N}$ denotes the number of distributions, with $M-1$ of them following a Gaussian distribution and one following a Dirac-delta function, as illustrated in the black box in Fig. \ref{fig6}. The transition probability matrix $A \in \mathbb{R}^{N_{\lambda} \times N_{\lambda}}$ specifies the probability of moving from each of the $N_{\lambda}$ hidden states to every other $N_{\lambda}$ hidden state, illustrated by the arrows among the hidden states in Fig. \ref{fig6}. The emission probability matrix $B \in \mathbb{R}^{N_{\lambda} \times M_{\lambda}}$ defines the probability that any of the $N_{\lambda}$ hidden states emit an observation. The observations correspond to the number of delayed samples caused by network constraints from the $M_{\lambda}$ distributions $(\mu,\sigma)$, shown as arrows from hidden states to the distributions in Fig. \ref{fig6}. Lastly, $\pi \in \mathbb{R}^{N_{\lambda}}$ indicates the probability of being in each hidden state at the initial time instant ($k=0$), represented by arrows from $\pi$ to all hidden states in Fig. \ref{fig6}.

Moreover, there is a finite set $S$ such that $\mathbf{card}(S)=N_{\lambda}$. The probability of transitioning from hidden network state $i$ to hidden network state $j$ adheres to the constraint $\Sigma_{j=1}^{N_{\lambda}} a_{ij}^{\lambda}=1, \ \forall i$, as depicted in Fig. \ref{fig6}. These individual transition probabilities across all hidden states collectively form the transition matrix $A$. The vector $\pi$, subject to the condition $\Sigma_{i=1}^{N_{\lambda}} \pi_{i}=1, \ \forall i$, represents the initial probability of starting in hidden state $i$ is such that $\pi_i = P(s_0 =i), (i \in S)$ where $s_k$ is the hidden state at time sample $k$ and $s_k \in S$.

In addition, there is a finite set $D$ with $\mathbf{card}(D)=M_{\lambda}$. The emission of any observation from distribution $j$ given any hidden state $i$ satisfies the constraint $\Sigma_{j=1}^{M_{\lambda}} b_{ij}^{\lambda}=1, \ \forall i$, as illustrated in Fig. \ref{fig6}. These emission probabilities are assembled to construct the emission probability matrix $B$. This matrix is expressed according to $P(d_k \mid s_k , s_{k-1} , ..., s_0)=P(d_k \mid s_k)$  where $d_k$ is the most probable distribution responsible for the observation at time sample $k$ and $d_k \in D$. 

In the SCHMM, $M_{\lambda}$ and $B$ correspond to a set of continuous distributions, unlike conventional HMMs where $M_{\lambda}$ and $B$ represent a discrete set of observations. While these continuous distributions increase the complexity of the training phase, they enable time delays to span a continuum of values and packet dropouts to be represented by a discrete value. The process of training the SCHMM and estimating the parameters that adjust the model to replicate the behavior of the actual communication network is described later in this work.

The sole purpose of this SCHMM is to provide $\hat{\tau}_{ij}(k)$, which is a prediction of the delay experienced by agent $i$ in information received from agent $j$ at the time instant $k$. This prediction will help agent $i$ to compute $\hat{x}_{ij}(k)$ using $x_{ij}(k-\tau_{ij})$ and $\hat{\tau}_{ij}(k)$ given that the dynamics and input commands of agent $j$, $\overrightarrow{u}_j$ are known to agent $i$. The input commands of agent $j$, $\overrightarrow{u}_j$ represent the latest optimized set of inputs computed from the SCHMM-LMPC at the time of sending the packet to agent $i$. Which then follows that each packet frame in this study sent from agent $j$ to agent $i$ includes $(x_{ij}(k-\tau_{ij}), \overrightarrow{u}_j)$. Subsequently, this will give rise to the following error term at agent $i$ for agent $j$:

\begin{equation} \label{eq1002}
E_{ij}(k)=x_j(k)-\hat{x}_j(k)
\end{equation}

\noindent given that $x_j(k)$ is unknown at agent $i$, and it has to depend on the accuracy of the predictions $\hat{x}_j(k)$ to drive the error term in \eqref{eq1002} to zero. Later discussions in this study leads to answering the question of the accuracy of the predictions made by agent $i$ about its neighboring agents. 

However, it has to be noted here that we now arrive at two error terms. The first error term is the error due to the distributed nature of the topology considered in this study shown in \eqref{eq13}. And the second error term is the error due to the network imperfections such as time delays and packet dropouts shown in \eqref{eq1002}. Using these two error terms, one can define the following:

\begin{equation} \label{eq1003}
E_{\lambda}(k)=E_i(k) + E_{ij}(k)
\end{equation}

\noindent where $E_{\lambda}(k)$ is the superposition of both error terms we arrived at, in this study. Now, the error term due to the network imperfections was tackled before in a previous study \cite{36} for the case of the single agent, not a MAS. In which, study \cite{36} used three algorithms to drive this error term to zero. The first algorithm is responsible for the initialization of the model $\lambda$. The second algorithm is the EM algorithm which is responsible for the tuning of the parameters of the model $\lambda$. The first and the second algorithms are offline, which means that these algorithms are executed pre-simulation. The third algorithm is the Viterbi algorithm which is responsible for the propagation of the predicted delays and this algorithm is online and executed during the simulation time.

Study \cite{36} addressed the error arising from network imperfections by learning parameters that best capture network behavior and employing an online algorithm to propagate predictions, effectively mitigating this error source. In contrast, the error introduced by the distributed topology, as examined in this study, is fundamentally different. Unlike network behavior, which is learnable offline, the topology varies across simulations and cannot be captured through static learning, provided it does not converge to an oscillating or alternating behavior. Consequently, enabling the SCHMM to also mitigate the topological error requires an EM algorithm capable of updating model parameters dynamically during simulation. The proposed incremental EM algorithm fulfills this need and represents the principal contribution of this work.

Building upon previous work \cite{36}, the initialization algorithm and the offline EM algorithm will be used in this study to arrive at $\lambda^*$, the offline optimized parameters of the SCHMM. This then shapes the algorithms presented in this study to be two algorithms. The first algorithm described in Algorithm \ref{alg1} is the incremental EM which tunes the parameters of the model $\lambda$ during the simulation and thus this algorithm is online. Moreover, the second algorithm shown in Algorithm \ref{alg2} is the Viterbi algorithm which is responsible for the propagation of the delay predictions and this algorithm is online. The Viterbi algorithm is adapted from \cite{36} and changed in such a way to apply for the MAS case, as it was previously used for a single control system only. These algorithms are described in detail later in this study. 

\section{Theoretical Analysis}
\label{sec3}
This study now considers the analysis for the presented methodology and divides it into two subsections. The first subsection focuses on the adaptation of the SCHMM and the introduction of the incremental EM algorithm. The second subsection presents the control design scheme using the SCHMM-LMPC and arrives at the theorem by which consensus of the MAS considered in this study is guaranteed.

\subsection{Adapted SCHMM using Incremental Expectation Maximization Algorithm}

To arrive at the incremental EM algorithm logically, we will first discuss the general method by which the parameters of the SCHMM are tuned to fit the behavior of the network.

The SCHMM plays a central role in this study, as it provides a probabilistic framework for modeling the stochastic behavior of communication networks subject to delays and packet dropouts. The SCHMM parameters capture the transition dynamics between hidden states as well as the probabilistic relationship between states and observed delays. However, in the context of distributed MASs, the communication topology itself introduces an additional source of uncertainty. Unlike network-induced imperfections, which can be learned from past data through offline training, the topological error is dynamic and changes with the graph structure in real time. This makes it essential to develop adaptive mechanisms that can continuously refine the SCHMM parameters during operation.

The model $\lambda$ parameters are categorized into two sets: the predefined parameters $(N_{\lambda}, M_{\lambda})$ and the trained parameters ($A, B, \pi$). The predefined ones are considered available beforehand, whereas the trained parameters are optimized by maximizing the log-likelihood function through the EM algorithm. Beginning with an initial $\lambda_0$, the EM iteratively updates the parameters to produce $\lambda'$, continuing until convergence to the optimal $\lambda^*$ within tolerance $\epsilon$. These final parameters are then employed to characterize the communication network dynamics by propagating the SCHMM using the Viterbi algorithm, which determines the most probable hidden state sequence given a set of observations. All while, executing the incremental EM algorithm to tune the parameters further adapting to the dynamics of the distributed topology.

The SCHMM is characterized by three principal problems that collectively define the system and capture network behavior. The first is the evaluation problem, which computes the probability that an observation sequence $O$ is generated by a model $\lambda$, expressed as $P(O \mid \lambda)$. This probability serves as the objective function and forms the basis for optimization, handled through both EM and incremental EM algorithms. The second is the decoding problem, which determines the most likely hidden state sequence explaining the observed data $O$, solved using the Viterbi algorithm. The third is the learning problem, concerned with updating the model parameters $(A, B, \pi)$ to maximize the objective function defined in the evaluation problem. This is achieved through EM and incremental EM, where optimization of the log-likelihood exploits the transformation of products into summations.

Now, that the general approach for calibrating the SCHMM parameters, to reflect network behavior, has been addressed. We start by recalling some of the key ideas in traditional EM algorithm that will help in developing the novel incremental EM algorithm for SCHMMs. As well as, highlight the main limitation that the EM algorithm suffers from which drives the development of the incremental EM algorithm in this study.

Offline EM algorithm outlines the computation of two main variables. The forward variable, $\alpha_t(i)$, denotes the probability of generating delay and dropout sequence $\{\tau_0, \dots, \tau_t\}$ while being in state $i$. Its recursive evaluation relies on the terms $a_{ji}^{\lambda}$ and $b_i^{\lambda}(\tau_t)$, where $a_{ji}^{\lambda}$ corresponds to the state transition probability and $b_i^{\lambda}(\tau_t)$ designates the observation probability in state $i$. Conversely, the backward variable, $\beta_t(i)$, expresses the probability of producing the remaining sequence $\{\tau_{T-1}, \dots, \tau_0\}$ given that the system occupies state $i$. It is assumed that dropouts in the set $\tau$ are encoded with a masking value $\Gamma=1/\epsilon$. Moreover, the variable $\zeta_t(i)$ refers to the posterior probability of occupying state $i$, while $\zeta_t(i,j)$ denotes the joint probability of transitioning from state $i$ to state $j$ at time $t+1$. Lastly, $\xi_t(i,g)$ captures the probability of selecting $g$-th Gaussian component from $(\mu,\sigma)$ in state $i$ for modeling the delay or dropout at time $t$. Further details are found in \cite{36} and references therein.

These variables allow the EM algorithm to tune the parameters of the SCHMM in order to arrive at $\lambda^*$. Although, it has to be noted that within distributed MASs, uncertainty also arises from the topology. While network-induced imperfections can be addressed through offline training on historical data, the topological error is time-varying and directly dependent on the instantaneous configuration of the graph, which requires a form of online learning. However, the method by which the EM algorithm tunes the parameters is inherently inapplicable for real-time scenarios. As it heavily relies on the computation of the forward and backward variables, and these are recursively calculated using the dataset of time delays and packet dropouts. This draws the attention to a limitation in the traditional EM algorithm, and dictates that it cannot solve the problem at hand on its own. Nevertheless, the EM algorithm cannot be ignored as it would be very computationally expensive to replace the algorithm with a new online one. Rather, this study will focus on building upon the EM algorithm and use it to train the SCHMM and arrive at $\lambda^*$, then develop an incremental setup of the EM algorithm in order to tackle the limitation.

Now, that the development of the incremental EM algorithm has been justified. The remaining part of this subsection is dedicated to Algorithm \ref{alg1} describing the purpose and details of the incremental EM algorithm. Discussion of Algorithm \ref{alg2}, which is the Viterbi algorithm, elaborates its complementary use in the delay predictions.


\begin{algorithm}[!t]
\begin{small}
\caption{Incremental Expectation-Maximization for SCHMM Parameter Update}
\KwIn{$\cdot$ Initial model parameters $\lambda_* = (\pi, A, B, \mu, \sigma)$\\
$\cdot$ learning rate $\eta$\\

$\cdot$ newly received packet of $x_j(k-\tau_{ij})$ and $\overrightarrow{u}_j$}
\KwOut{Updated model parameters $\lambda^*$}

\Begin{
$\tau_{\text{prev}} = \arg\min_{\tau} \left\| x_j(k-\tau_{ij}) - \hat{x}_j(k-\tau \mid \overrightarrow{u}_j) \right\|$

\BlankLine
\For{each incoming observation computation $\tau_{\text{prev}}$}{
    
    \For{each state $i$}{
        $\gamma_t(i) \leftarrow \frac{\pi_i b_i^{\lambda}(\tau_{\text{prev}})}{\sum\limits_{j=1}^{N_{\lambda}} \pi_j b_j^{\lambda}(\tau_{\text{prev}})}$
    }

    \For{each state pair $(i,j)$}{
        $\gamma_t(i,j) \leftarrow \frac{\pi_i a_{ij}^{\lambda} b_j^{\lambda}(\tau_{\text{prev}})}{\sum\limits_{p=1}^{N_{\lambda}} \sum\limits_{q=1}^{N_{\lambda}} \pi_p a_{pq}^{\lambda} b_q^{\lambda}(\tau_{\text{prev}})}$
    }
    
    \BlankLine
    
    \For{each state $i$}{
        $\pi_i \leftarrow (1 - \eta)\pi_i + \eta\, \gamma_t(i)$
    }
    
    \For{each state pair $(i,j)$}{
        $a_{ij}^{\lambda} \leftarrow (1 - \eta)a_{ij}^{\lambda} + \eta\, \gamma_t(i,j)$
    }
    
    \For{each state $i$ and distribution $D_i$}{
        $\mu_{i} \leftarrow (1 - \eta)\mu_{i} + \eta\, \gamma_t(i) \tau_{\text{prev}}$ \\
        $\sigma_{i}^2 \leftarrow (1 - \eta)\sigma_{i}^2 + \eta\, \gamma_t(i) (\tau_{\text{prev}} - \mu_{i})^2$
    }
    
    Update $\lambda_{\text{old}}^* \leftarrow \lambda_{\text{new}}^*\;$
}}
\label{alg1}
\end{small}
\end{algorithm}

It has to be noted that the purpose of the incremental EM algorithm is to ensure that the dynamic error due to the topology is mitigated. The scheme by which the SCHMM suppresses the error due to the network imperfections is the same for the error due to the distributed topology. The only difference is that the behavior of the network can be learned offline using the traditional EM algorithm. But, for the error due to the distributed topology, the graphs are dynamical and native to the time instant of the simulation. This was the driving factor of this study to introduce the incremental form of the EM algorithm which allows the online learning of the parameters of the SCHMM which in turn suppresses the dynamical error due to the distributed topology. Given that the behavior of the error shown in \eqref{eq1003} is not an oscillating or alternating one which is valid in our case as a spanning tree-undirected graph is assumed. So as the incremental EM algorithm can achieve convergence in a finite time and not to chase an ever-changing target.

Algorithm \ref{alg1} implements an incremental EM procedure for the SCHMM, designed to refine the model parameters in real time as new packet information becomes available. Unlike the offline EM, which processes the entire dataset at once, the incremental version updates the parameters $(\pi, A, B, \mu, \sigma)$ continuously. This online adaptation enables the SCHMM to capture dynamic variations in the distributed topology without reprocessing past observations. The algorithm balances previously learned information with new evidence using a learning rate $\eta$, ensuring responsiveness. In doing so, it mitigates error terms arising from network and topology dynamics that cannot be addressed by offline training alone.


The algorithm begins by estimating the effective delay $\tau_{\text{prev}}$, which measures the discrepancy between the actual received delayed state and its predicted counterpart, thereby providing the observation that drives parameter updates. For each new observation, the E-step computes the posterior state probabilities $\gamma_t(i)$, indicating how likely it is that the system is in state $i$, and the transition probabilities $\xi_t(i,j)$, estimating the likelihood of moving from state $i$ to $j$, with normalization across all states. In the M-step, the initial distribution $\pi_i$ is incrementally updated with $\gamma_t(i)$ to capture the evolving likelihood of starting in each state, while the transition probabilities $a_{ij}^{\lambda}$ are refined using $\xi_t(i,j)$ to adjust the model’s representation of state-to-state dynamics under current network conditions. For each state, the emission distribution parameters $(\mu_i, \sigma_i)$ are updated with the latest delay observation weighted by $\gamma_t(i)$, allowing the model to adapt the mean and variance of delays to observed data. Finally, the model parameters are refreshed to yield $\lambda^*_{\text{new}}$, ensuring that the SCHMM evolves continuously and accurately reflects real-time communication delays and packet dropouts in a distributed system.

One remark here is that Algorithm \ref{alg1} appears to be computationally expensive. However, as discussed in \cite{21}, in order to best mimic the behavior of the network, it is sufficient to have few hidden states and distributions in the SCHMM. As increased hidden states leads to granulation and over fitting. This then follows that Algorithm \ref{alg1} would not be computationally-consuming.

\begin{algorithm}[!t]
\begin{small}
\caption{Adapted Viterbi Algorithm}
\KwIn{$\cdot$ Optimized SCHMM parameters $\lambda^* = (\pi^*, A^*, B^*)$\\

$\cdot$ newly received packet of $x_j(k-\tau_{ij})$ and $\overrightarrow{u}_j$}
\KwOut{Predicted delay $\tau_{\text{next}}$}

\Begin{
$\tau_{\text{prev}} = \arg\min_{\tau} \left\| x_j(k-\tau_{ij}) - \hat{x}_j(k-\tau \mid \overrightarrow{u}_j) \right\|$

\textbf{Step 1: Update the most likely state sequence with the previous delay} \\
\Begin{
    Using the newly received packet delay $\tau_{\text{prev}}$:
    \[
    s_{k-1} = \arg \max_{s_k} \alpha_{k-1}(s_{k-1}) \beta_{k-1}(s_{k-1})
    \]
}

\textbf{Step 2: Predict the next state $s_{k}$} \\
\Begin{
    Using the updated state $s_{k-1}$ and transition probabilities $A^*$:
    \[
    s_{k} = \arg \max_j a^{\lambda *}_{s_{k-1},j}
    \]
}

\textbf{Step 3: Predict the next delay $\tau_{\text{next}}$} \\
\Begin{
    Using the predicted state $s_{k}$ and the Gaussian and Dirac-delta distributions from $B^*$ and $(\mu, \sigma)$:
    \[
    \tau_{\text{next}} = \arg \max_{\tau} b^{\lambda *}_{s_{k}}(\tau)
    \]
}

\Return $\tau_{\text{next}}\;$}
\label{alg2}
\end{small}
\end{algorithm}

Algorithm \ref{alg2} presents an adapted version of the Viterbi algorithm for delay prediction in distributed networks. Its primary role is to exploit the optimized SCHMM parameters $(\pi^*, A^*, B^*)$ to estimate the most likely hidden state sequence and forecast the next communication delay $\tau_{\text{next}}$. By combining the state transition dynamics with the probabilistic emission distributions, the algorithm ensures that predictions are consistent with both past observations and the learned model. This adaptation makes the Viterbi algorithm suitable for online operation in multi-agent systems, where communication imperfections such as delays and dropouts are time-varying and uncertain.

The algorithm starts by estimating the most recent delay $\tau_{\text{prev}}$ through comparison between the received delayed state and its predicted counterpart. In the first step, the most probable previous state sequence is updated using the forward and backward variables $\alpha_{k-1}(s_{k-1})$ and $\beta_{k-1}(s_{k-1})$. Next, the subsequent state $s_k$ is predicted by selecting the state that maximizes the transition probability $a^{\lambda *}_{s_{k-1},j}$. Finally, the upcoming delay $\tau_{\text{next}}$ is obtained by evaluating the emission distribution of the predicted state $s_k$, selecting the delay value with the highest likelihood under $b^{\lambda *}_{s_k}(\tau)$. By iterating this process with each new packet, the algorithm provides accurate, real-time delay predictions that support consensus in distributed control systems.

It is essential to note that the work of this study extends and improves the previous work presented in \cite{36} where the focus was on a single control system and the delay experienced was in the feedback and forward channels. This study deals with two new traits that were missing from \cite{36}, that are now alleviated in this new presented work. The first change was the introduction of the incremental EM algorithm which allows the online tuning of the parameters of the model $\lambda$ which was missing from the work in the literature and is regarded as the main contribution of this study. The second change is the adaptation of the Algorithm \ref{alg2} to the MAS case, as the earlier presentations were limited to solve the case of a single agent or a single control system.

\subsection{LMPC}

Now, the LMPC design is considered where, a theorem is presented to achieve the design of the matrix $K$ which will result in the consensus of the MAS and the minimization of the compact cost function for each agent $i$ using their version of the global problem.

\textbf{Theorem 1}: \textit{The error in the controller system \eqref{eq16} is asymptotically stable if there exists a constant $\alpha>0$, a matrix $\Pi \in \mathbb{R}^{mN \times nN}$ and a positive definite symmetric matrix $\Omega \in \mathbb{R}^{Nn \times nN}$ such that:}

\begin{align*} 
\min_{\text{subject to:}}\alpha
\end{align*}

\begin{equation} \label{eq220}
\begin{bmatrix}
1 & E_i^T(k)\\ 
E_i(k) & \alpha^{-1}P_v
\end{bmatrix}\geq0
\end{equation}

\begin{equation} \label{eq221}
\begin{scriptsize}
\begin{bmatrix}
(\Omega-I)\Omega & \Omega(A_e^{-1})^TA_c^TA_e^T+\Omega(\Omega^{-1})^T\Pi^TB_c^TA_e^T\\ 
A_eA_cA_e^{-1}\Omega+A_eB_c\Pi & -\Omega^{-1}
\end{bmatrix}<0
\end{scriptsize}
\end{equation}

\noindent \textit{this follows that the design matrix $K$ in \eqref{eq21} is:}

\begin{equation} \label{eq26}
K=\Pi\Omega^{-1}
\end{equation}

\textbf{Proof of Theorem 1:} Assume the following Lyapunov candidate function as:

\begin{equation} \label{eq22}
V(E_i(k))=E_i^T(k)P_vE_i(k)
\end{equation}

\noindent with $P_v$ is being a real symmetric positive definite matrix.

We also consider that:

\begin{equation} \label{eq23}
\begin{gathered}
\Delta V(E_i(k))=V(E_i(k+p+1))-V(E_i(k+p)) \leq -\\ (E_i^T(k+p)P_JE_i(k+p)+U_{i}^T(k+p)Q_JU_{i}(k+p))
\end{gathered}
\end{equation}

\noindent and summing on $p$ from $0$ to $\infty$, we get:

\begin{equation} \label{eq24}
\begin{gathered}
\sum_{p=0}^{\infty}V(E_i(k+p+1))-V(E_i(k+p)) \leq -\\ \sum_{p=0}^{\infty} (E_i^T(k+p)P_JE_i(k+p)+U_{i}^T(k+p)Q_JU_{i}(k+p))
\end{gathered}
\end{equation}

Right hand side of inequality \eqref{eq24} is $J_i(k)$ while summing to $\infty$ imposes that $E_i(k+\infty)\rightarrow0$ due to the asymptotic convergence. Thus, one gets:

\begin{equation} \label{eq25}
-V(E_i(k)) \leq -J_i(k) \rightarrow J_i(k) \leq V(E_i(k))
\end{equation}

Thus, the optimization problem reduces to:

\begin{equation} \label{eq26}
\min_{U_{i}(k)} V(E_i(k))
\end{equation}

\noindent we also then assign an upper bound so that:

\begin{equation} \label{eq27}
V(E_i(k)) \leq \alpha
\end{equation}

And, the optimization problem becomes:

\begin{equation} \label{eq28}
\min_{V(E_i(k))} \alpha
\end{equation}

We then consider the derivation of LMIs \eqref{eq220} and \eqref{eq221}, starting by substituting \eqref{eq22} in \eqref{eq27} to get:

\begin{equation} \label{eq29}
\begin{gathered}
E_i^T(k)P_vE_i(k) \leq \alpha\\
\alpha^{-1}E_i^T(k)P_vE_i(k) \leq 1\\
1-\alpha^{-1}E_i^T(k)P_vE_i(k) \geq 0
\end{gathered}
\end{equation}

We now apply Schur's complement to get:

\begin{equation} \label{eq30}
\begin{bmatrix}
        1 & E_i^T(k)\\ 
        E_i(k) & \alpha^{-1}P_v
    \end{bmatrix}\geq 0
\end{equation}

\noindent this defines LMI \eqref{eq220}.

Now, we consider \eqref{eq16} and \eqref{eq21} to get:

\begin{equation} \label{eq31}
\begin{gathered}
E_i(k+p+1)=A_eA_cA_e^{-1}E_i(k+p)+A_eB_cKE_i(k+p)\\
E_i(k+p+1)=(A_eA_cA_e^{-1}+A_eB_cK)E_i(k+p)
\end{gathered}
\end{equation}

Then, we should focus on $\Delta V(E_i(k+p))$ and using \eqref{eq31}, we get:

\begin{equation} \label{eq32}
\begin{gathered}
\Delta V(E_i(k+p))=V(E_i(k+p+1))- V(E_i(k+p)) \\
=E_i^T(k+p+1)P_vE_i(k+p+1)-E_i^T(k+p)P_vE_i(k+p)\\
=E_i^T(k+p)(A_eA_cA_e^{-1}+A_eB_cK)^TP_v(A_eA_cA_e^{-1}+\\
A_eB_cK)E_i(k+p)-E_i^T(k+p)P_vE_i(k+p)\\
=E_i^T(k+p)((A_eA_cA_e^{-1}+A_eB_cK)^TP_v(A_eA_cA_e^{-1}+\\
A_eB_cK)-P_v)E_i(k+p)
\end{gathered}
\end{equation}

\noindent and then ensuring that $\Delta V(E_i(k+p)) < 0$, one gets:

\begin{equation} \label{eq33}
\begin{gathered}
(A_eA_cA_e^{-1}+A_eB_cK)^TP_v(A_eA_cA_e^{-1}+\\
A_eB_cK)-P_v<-I<0    
\end{gathered}
\end{equation}

We now apply Schur's complement to get:

\begin{equation} \label{eq34}
\begin{bmatrix}
        I-P_v & (A_eA_cA_e^{-1}+A_eB_cK)^T\\ 
        * & -P_v
    \end{bmatrix}< 0
\end{equation}

\noindent now, it is wise to pre and post multiply inequality \eqref{eq34} by $\text{diag}\{P_v^{-1},I\}>0$ to get:

\begin{equation} \label{eq35}
\begin{scriptsize}
\begin{bmatrix}
        (P_v^{-1}-I)P_v^{-1} & P_v^{-1}(A_eA_cA_e^{-1}+A_eB_cK)^T\\ 
        (A_eA_cA_e^{-1}+A_eB_cK)P_v^{-1} & -P_v
    \end{bmatrix}< 0
\end{scriptsize}
\end{equation}

Then, it will be helpful to define the following:

\begin{equation} \label{eq36}
\Omega \triangleq P_v^{-1} ~ \text{and} ~ \Pi \triangleq K\Omega
\end{equation}

\noindent and using these definitions we arrive at LMI \eqref{eq221}:

\begin{equation} \label{eq37}
\begin{scriptsize}
\begin{bmatrix}
        (\Omega^{-1}-I)\Omega & \Omega(A_e^{-1})^TA_c^TA_e^T+\Omega(\Omega^{-1})^T\Pi^TB_c^TA_e^T\\
        A_eA_cA_e^{-1}\Omega+A_eB_c\Pi & -\Omega^{-1}
    \end{bmatrix}< 0
\end{scriptsize}
\end{equation}

\noindent which then concludes that design matrix $K$ is:

\begin{equation} \label{eq38}
K=\Pi \Omega^{-1}
\end{equation}

\noindent proving that the error in controller system \eqref{eq16} is asymptotically stable which then completes the proof. $~~~~~~~~~~~~~~~~~~~~~~~~~\blacksquare$

Trace-backing through the equations will allow each agent $i$ to arrive at $\overrightarrow{u}_i$ which is the optimized input commands for the prediction window. This then drives the agent $i$ towards $\delta_i(k)$ and eventually converge to the state of consensus for all agents.

\section{Numerical Simulations}
\label{sec4}
Two examples are simulated to validate the theoretical results and demonstrate the practical applicability of the proposed SCHMM-LMPC. These include one example where agents are a part of a centralized topology which goes to show the versatility of the proposed approach that is not only native to the specific problem it was tailored to solve, but also, to other setups and topologies as well. The second example presents the case where the agents are a part of a distributed topology with no centralization in the graph which goes to showcase the capability of the proposed SCHMM-LMPC approach in solving the distributed optimal consensus problem in a MAS under network imperfections, solving all three problems discussed earlier in the literature review.

The network imperfections, with both its significant factors, time delays and packet dropouts, are simulated first using NS-2 and then the SCHMM is trained using the offline EM algorithm. After that, the two consensus simulations are performed once in a centralized topology and once in a distributed topology, while performing the online algorithms (Algorithm \ref{alg1} and Algorithm \ref{alg2}). This then dictates that this section is divided into the SCHMM training subsection and the SCHMM-LMPC consensus simulations with two examples subsection.

\subsection{SCHMM Training}

Since the foundation of the proposed approach in this study lies in the SCHMM, it is appropriate to begin by addressing the dataset used to train it. The initial stage involves conducting network simulations in NS-2 to gather raw measurements of time delays and packet dropouts. These measurements constitute the training dataset for the SCHMM, enabling it to adapt and replicate the characteristics of the simulated network environment. The dataset was collected from a typical network used in a distributed MAS (ad-hoc communication) where the agents were performing maneuvers past each other to finally reach a consensus point. Following this, parameter initialization is performed. For the SCHMM, with $N_{\lambda}=3$ and $M_{\lambda}=4$, and matrices ($A, B, \pi$) initialized uniformly, K-means clustering is applied to determine the initial means $\mu$ and variances $\sigma$ of the distributions.

And upon executing the offline EM algorithm, with a maximum of $50$ iterations and a convergence threshold $\epsilon=10^{-8}$, the optimized parameters $\lambda^*$ were obtained as follows:

\begin{small}

\vspace{1mm}

$\pi^*=\begin{bmatrix}
0.4215 & 0.4572 & 0.1213
\end{bmatrix}$

\vspace{1mm}

$A^*=\begin{bmatrix}
0.6832 & 0.2079 & 0.1089\\
0.2894 & 0.5538 & 0.1568\\
0.1245 & 0.3761 & 0.4994 
\end{bmatrix}$

\vspace{1mm}

$B^*=\begin{bmatrix}
0.0221 & 0.4528 & 0.3647 & 0.1604\\
0.0213 & 0.5327 & 0.2934 & 0.1526\\
0.0198 & 0.5021 & 0.3504 & 0.1277
\end{bmatrix}$

\vspace{1mm}

$\mu^*=\begin{bmatrix}
46.00 & 49.85 & 58.17 & 100000
\end{bmatrix}$

\vspace{1mm}

$\sigma^*=\begin{bmatrix}
0.4149 & 1.0733 & 2.9872 & 0.0001
\end{bmatrix}$

\vspace{1mm}

\end{small}

\begin{figure}[!t]
\centering
\includegraphics[width=0.17\textwidth,height=0.13\textwidth]{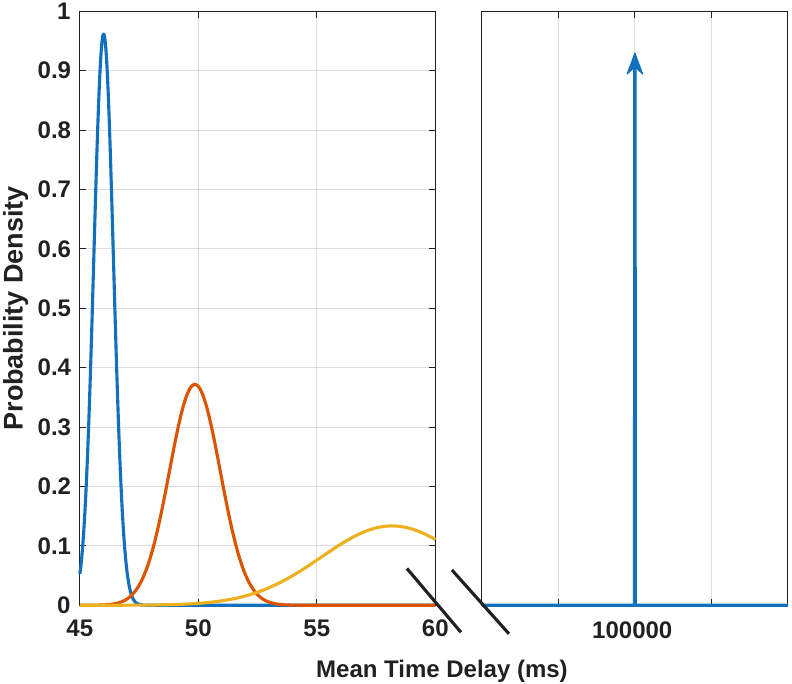}
\caption{Gaussian distributions corresponding to time delays measured in ms and Dirac-delta function corresponding to packet dropouts masked with the value $10^5$ms.}
\label{fignetw}
\end{figure}

The outcomes of applying the offline EM algorithm to the SCHMM are depicted in Fig. \ref{fignetw}, where the three Gaussian distributions correspond to time delays, while the Dirac-delta function represents packet dropouts. The mean of the Dirac-delta distribution corresponds to the masking value assigned to packet dropouts in the training set, namely $10^5$ms. As expected, its variance is negligibly small, reflecting the discrete nature of packet dropouts. These model parameters are then copied to all agents in the simulations, and upon running Algorithm \ref{alg1} $(\eta = 0.1)$, model's parameters converge to different values in different agents as the model reflects the behavior of the network and the topology locally to each agent. This calibrated model was subsequently employed to generate new instances of varying time delays and packet dropouts using Algorithm \ref{alg2}.

The effect of applying Algorithm \ref{alg1} can be further demonstrated by analyzing the SCHMM parameters of agent $3$ in Example $2$ which is discussed below. After $10$ seconds of simulation, the resulting parameters reveal that the distributed nature of the system induces notable variations. These variations are significant, as evidenced by both the number of parameters affected and the magnitude of their percentage changes. Such results emphasize the relevance of the proposed approach, which provides a systematic means to mitigate errors arising from the distributed topology as well as from network-induced imperfections.

\begin{small}

\vspace{1mm}

$\pi^*=\begin{bmatrix}
0.4585 & 0.4452 & 0.0963
\end{bmatrix}$

\vspace{1mm}

$A^*=\begin{bmatrix}
0.6735 & 0.2170 & 0.1095\\
0.2950 & 0.5477 & 0.1573\\
0.1301 & 0.3704 & 0.4995 
\end{bmatrix}$

\vspace{1mm}

$B^*=\begin{bmatrix}
0.0250 & 0.4480 & 0.3680 & 0.1590\\
0.0251 & 0.5279 & 0.2973 & 0.1497\\
0.0198 & 0.5111 & 0.3415 & 0.1276
\end{bmatrix}$

\vspace{1mm}

$\mu^*=\begin{bmatrix}
46.66 & 48.29 & 55.87 & 100000
\end{bmatrix}$

\vspace{1mm}

$\sigma^*=\begin{bmatrix}
0.4094 & 1.4031 & 2.8222 & 0.0001
\end{bmatrix}$

\vspace{1mm}

\end{small}

\subsection{SCHMM-LMPC}

Now, that the SCHMMs have been trained using the offline EM algorithm, the MAS simulation targeting consensus can be carried out. This will take place through two examples, the first describing a centralized MAS which is presented to show that the framework presented in this paper is not native to the case it solves, but rather a framework that can handle various topologies in MASs. The second example discussing a distributed topology of a MAS and this example is featured in this paper to act as a validation that the proposed methodology indeed solves the three problems discussed in the literature review.

\begin{figure}[!t]
\begin{minipage}[b]{0.4\linewidth}
\centering
\includegraphics[width=0.55\textwidth]{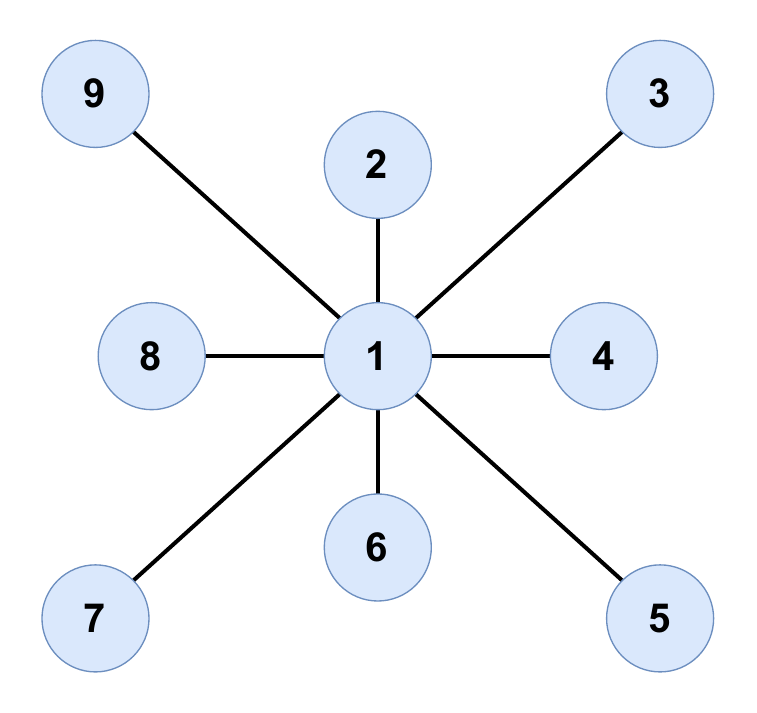}
\caption{Graph describing MAS in Example 1.}
\label{q1}
\end{minipage}
\hspace{0.5cm}
\begin{minipage}[b]{0.4\linewidth}
\centering
\includegraphics[width=\textwidth]{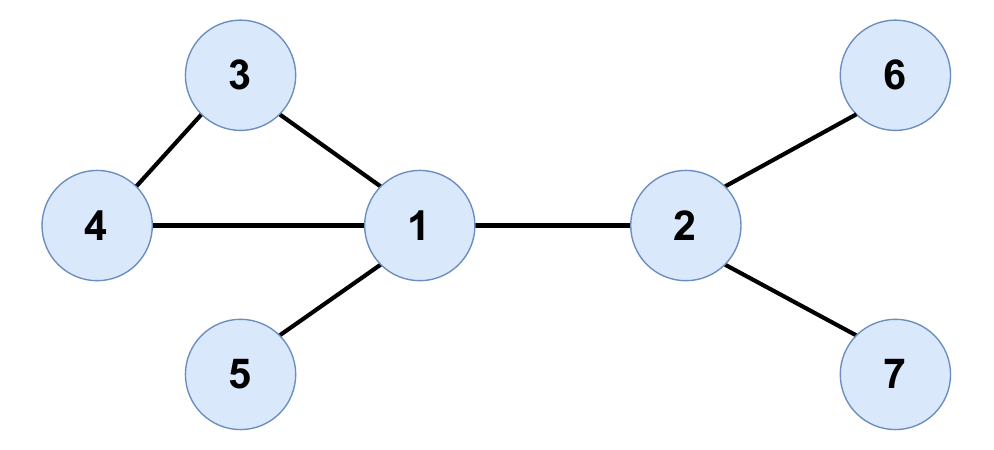}
\caption{Graph describing MAS in Example 2.}
\label{q2}
\end{minipage}
\end{figure}

\textbf{Example 1:} For the first numerical example, a centralized MAS is considered, where each agent's dynamics are modeled by \eqref{eq1}, $(n=6, m=4)$ with nine agents, $(N=9)$, as shown in fig. \ref{q1} depicting the graph used in the simulation.

Now, we define parameters such as $P_v =I_{nN}$ and $\alpha=10^{-3}$. Initial guesses are used for the $\Omega = I_{nN}+10^{-6}*I_{nN}$ and $\Pi=1_{mN \times nN}$. Random initial states are followed using randn(.) function in MATLAB. The SCHMM trained earlier was used to predict the network behavior during the simulation while running the online algorithms.

Two points to be made from the results of Example 1. The first is that the MAS achieved consensus quickly and smoothly before the $5$ second mark. This was the anchor example and the harder problem is presented in Example 2. The second, more significant, result is that the parameters of the model $\lambda$ did not change for the entire time of the simulation even while running Algorithm \ref{alg1}. This goes to prove that Algorithm \ref{alg1} is indeed responsible only for error due to the distributed topology, and has no effect when the topology is centralized. This goes to show the versatility of the approach and that it is not native to the problem at hand only, but rather a complete framework engulfing centralized and distributed topologies.

\textbf{Example 2:} For the second numerical example, a distributed MAS is considered, where same agent's dynamics used in Example 1 are used in this Example, however with seven agents, $(N=7)$, as shown in fig. \ref{q2} depicting the graph used in the simulation.

\begin{figure}[!t]
\centering
\begin{minipage}[b]{0.4\linewidth}
\centering
\includegraphics[width=\textwidth]{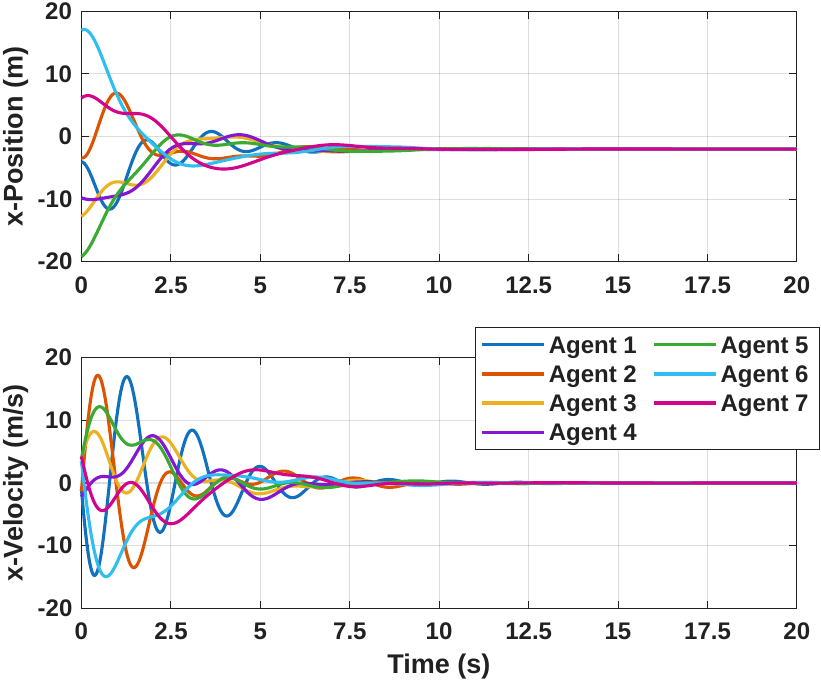}
\caption{State trajectories in the $x$-direction for Example 2.}
\label{num6}
\end{minipage}
\hspace{0.5cm}
\begin{minipage}[b]{0.4\linewidth}
\centering
\includegraphics[width=\textwidth]{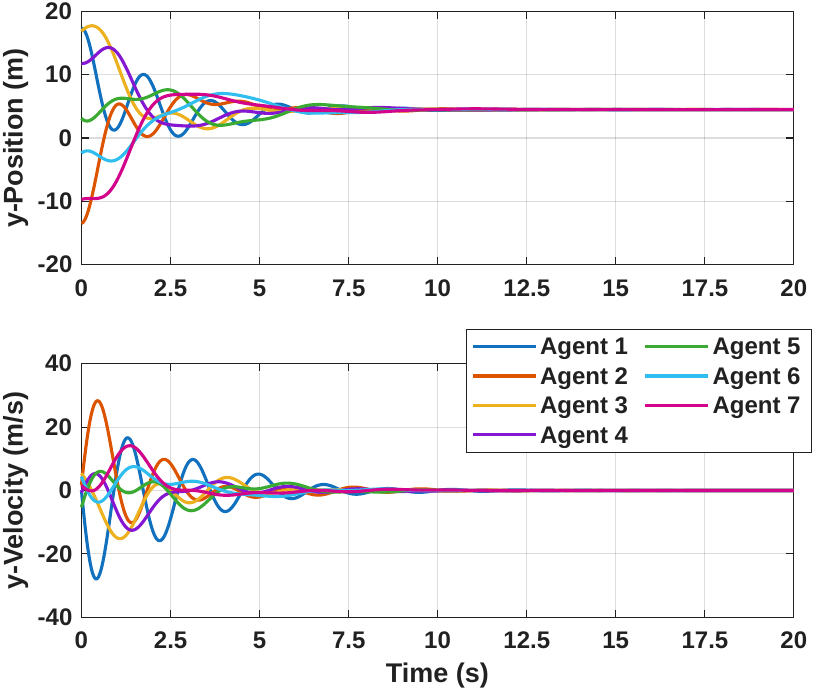}
\caption{State trajectories in the $y$-direction for Example 2.}
\label{num7}
\end{minipage}
\end{figure}

\begin{figure}[!t]
\centering
\begin{minipage}[b]{0.4\linewidth}
\centering
\includegraphics[width=\textwidth]{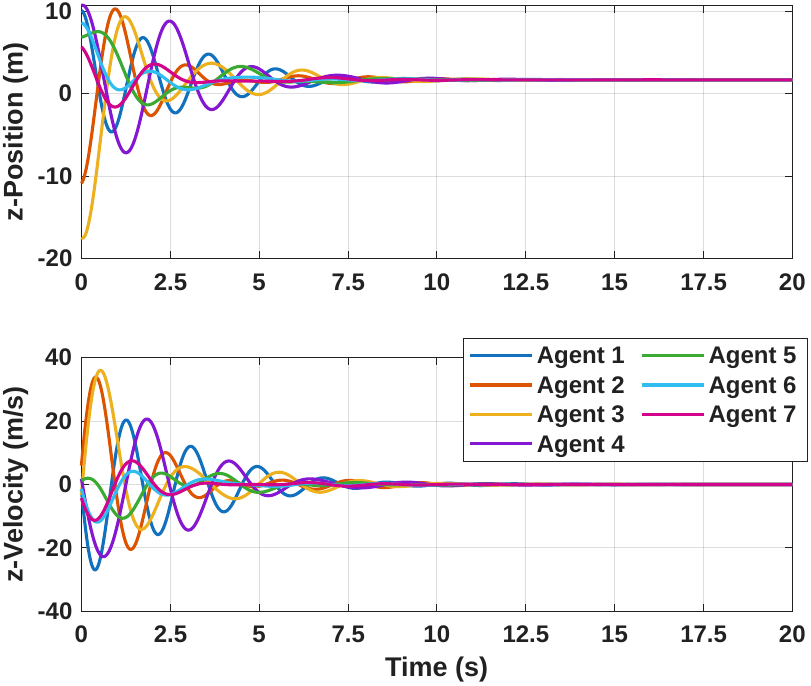}
\caption{State trajectories in the $z$-direction for Example 2.}
\label{num8}
\end{minipage}
\hspace{0.5cm}
\begin{minipage}[b]{0.4\linewidth}
\centering
\includegraphics[width=\textwidth]{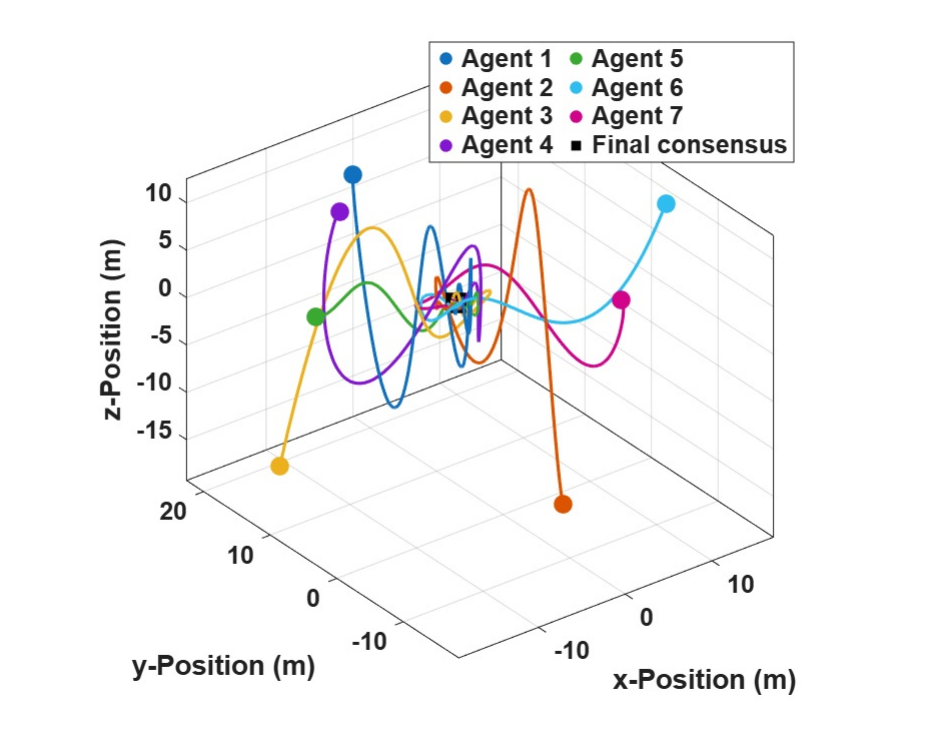}
\caption{State trajectories in the 3-D space for Example 2.}
\label{num9}
\end{minipage}
\end{figure}

\begin{figure}[!t]
\centering
\includegraphics[scale=0.2]{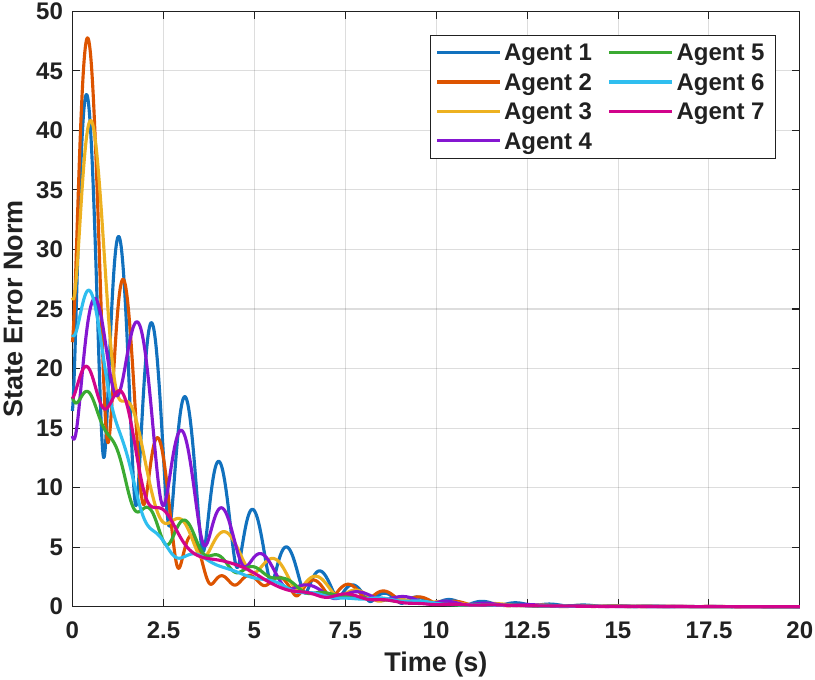}
\caption{State Error Norm between the position of each agent and its group term $\delta_i(k)$ against time in Example 2.}
\label{num10}
\end{figure}

We use the same parameters used in Example 1 as to be fair in the comparison between the change in the graph topology. The SCHMM trained earlier was used again to predict the network behavior during the simulation while running again the online algorithms.

The Figs. \ref{num6}, \ref{num7} and \ref{num8} show the results of the simulation of the distributed MAS where the evolution of the agents positions and velocities for $20$ seconds in the $x$, $y$ and $z$ directions, respectively. Fig. \ref{num9} shows the evolution of the agents in the 3-D space achieving consensus. Fig. \ref{num10} shows the state error norm for each agent in the simulation depicting its decay to zero indicating achieving consensus. The state error norm is defined as the 2-norm between the position of each agent with its local consensus point $\delta_i(k)$ for each agent $i$ in the MAS.

\section{Conclusion}
\label{sec5}
This paper presented the Semi-Continuous Hidden Markov Model – Lyapunov-based Model Predictive Control (SCHMM-LMPC) framework for distributed optimal control of Multi-Agent Systems (MASs) subject to time delays and packet dropouts. The integration of Lyapunov-based MPC with Linear Matrix Inequalities (LMIs) ensures optimality and consensus, while the SCHMM provides accurate delay prediction and compensation. The new formulation presented in this study defines two errors, one due to the network imperfections and the other due to lack of centralization. The main contribution is the incremental Expectation Maximization (EM) algorithm, which enables the SCHMM to update its parameters online and adapt to varying network conditions as well as adverse effects of the distributed topology. Unlike offline methods, this incremental EM mitigates both topology-induced errors and uncertainties from communication imperfections in real time. Numerical examples confirm the framework’s effectiveness in maintaining consensus under delays and dropouts, demonstrating its scalability and applicability to diverse MAS topologies. Overall, the study contributes a novel online adaptation mechanism that strengthens the link between predictive control and probabilistic modeling in distributed networks.

\bibliographystyle{IEEEtran}
\bibliography{example.bib}

\vspace{-2cm}

\begin{IEEEbiography}[{\includegraphics[width=1in,height=1.25in,clip,keepaspectratio]{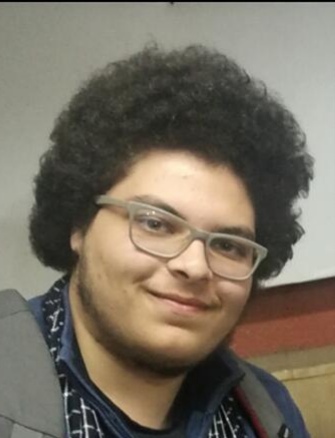}}]{Loaie Solyman} was born in Cairo, Egypt on June 1, 1998. He received a B.Sc. and M.Sc. degrees in mechatronics engineering in 2021 and 2024 from the German University in Cairo (GUC). 

He joined the mechatronics engineering department in the GUC as a teaching assistant from 2021 where he is currently a PhD candidate. His research interests include networked control systems and data-driven control theory.
\end{IEEEbiography}

\vspace{-2cm}

\begin{IEEEbiography}[{\includegraphics[width=1in,height=1.25in,clip,keepaspectratio]{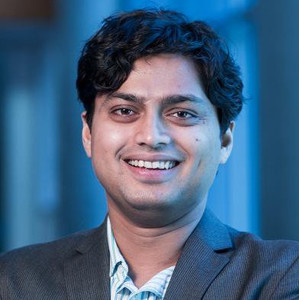}}]{Aamir Ahmad} obtained a B-Tech (with honors) degree in 2008 in Civil Engineering from the Indian Institute of Technology (IIT), Kharagpur, India. Ahmad obtained a PhD (with merit) degree in 2013 in electrical and computer engineering from the University of Lisbon, Portugal.

He is a tenure-track professor of Flight Robotics and the Deputy Director (Research) at the Institute for Flight Mechanics and Controls, Faculty of aerospace engineering and geodesy, University of Stuttgart, Germany. He is also a Research Group Leader at the Max Planck Institute for Intelligent Systems in Tübingen, where he was previously a research scientist (2016-2020). His research interests mainly include aerial robotics and multi-robot systems.
\end{IEEEbiography}

\vspace{-2cm}

\begin{IEEEbiography}[{\includegraphics[width=1in,height=1.25in,clip,keepaspectratio]{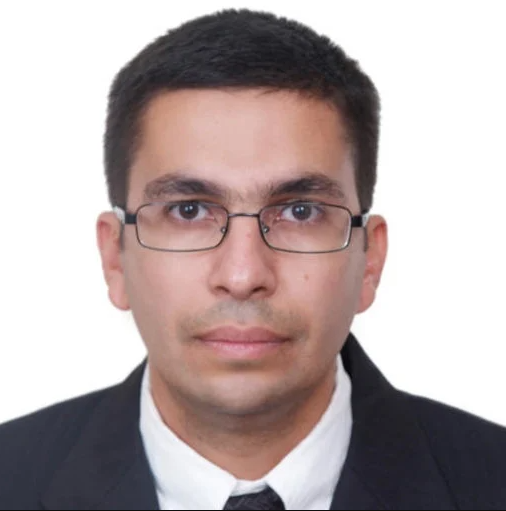}}]{Ayman El-Badawy} was born in Cairo, Egypt on August 17, 1972. He received a B.Sc. and M.Sc. degrees in mechanical engineering in 1993 and 1995 from the American University in Cairo. He received his Ph.D. in mechanical engineering in 2000 from Virginia Polytechnic Institute and State University. 

He has industrial experience with A.O. Smith Corporation, USA and Nile Aster, Egypt. He is currently the head of mechatronics engineering department at the German University in Cairo (GUC). His research interests include control theory and reinforcement learning.
\end{IEEEbiography}

\end{document}